\documentclass[12pt]{article}
\usepackage{amsmath}
\usepackage{amsfonts}
\usepackage{amssymb}
\usepackage[mathscr]{eucal}

\newif{\ifcomentarios}
\comentariosfalse

\newtheorem{theorem}{Theorem}

\newtheorem{lemma}[theorem]{Lemma}

\newtheorem{remark}[theorem]{Remark}

\newcommand{\be}{\begin{eqnarray}}
\newcommand{\en}{\end{eqnarray}}
\newcommand{\bee}{\begin{eqnarray*}}
\newcommand{\ene}{\end{eqnarray*}}

\input{tcilatex}

\begin{document}

\author{S. L. Carvalho\thanks{%
Work supported by FAPESP under grant \#06/60711-4. Present adrress:
Instituto de Ci\^{e}ncias e Tecnologia, UNIFESP, 11231-280 S\~{a}o Jos\'{e}
dos Campos, SP, Brasil. Email: \texttt{slcarvalho@unifesp.br}}, \ D. H. U.
Marchetti\thanks{%
Present Address: Mathematics Department, \ The University of British
Columbia, Vancouver, \ BC, \ Canada \ V6T 1Z2. Email: \texttt{%
marchett@math.ubc.ca}} \ \& W. F. Wreszinski\thanks{%
Email:\texttt{\ wreszins@fma.if.usp.br}} \\
Instituto de F\'{\i}sica\\
Universidade de S\~{a}o Paulo \\
Caixa Postal 66318\\
05314-970 S\~{a}o Paulo, SP, Brasil}
\title{On the uniform distribution of the Pr\"{u}fer angles and its
implication for a sharp spectral transition of Jacobi matrices with randomly
sparse perturbations}
\date{}
\maketitle

\begin{abstract}
In the present work we consider Jacobi matrices with random uncertainty in
the position of sparse perturbations. We prove (Theorem \ref{thethe}) that
the sequence of Pr\"{u}fer angles $(\theta _{k}^{\omega })_{k\geq 1}$ is u.d
mod $\pi $ for all $\varphi \in \left[ 0,\pi \right] $ with exception of a
set of rational multiples of $\pi $ and for almost every $\omega $ with
respect to the product $\nu =\prod_{j\geq 1}\nu _{j}$ of uniform measures on 
$\left\{ -j,\ldots ,j\right\} $. Together with an improved criterion for
pure point spectrum (Lemma \ref{l2}), this provides a simple and natural
alternative proof of a result of Zlatos (J. Funct. Anal. \textbf{207},
216-252 (2004)): the existence of pure point (p.p) spectrum and singular
continuous (s.c.) spectra on sets complementary to one another with respect
to the essential spectrum $\left[-2,2 \right]$, outside sets $A_{sc}$ and $%
A_{pp}$, respectively, both of zero Lebesgue measure (Theorem \ref%
{tmarsilwre}). Our method allows for an explicit characterization of $A_{pp}$%
, which is seen to be also of dense p.p. type, and thus the spectrum is
proved to be exclusively pure point on one subset of the essential spectrum.

\medskip

\end{abstract}

\section{Introduction}

\setcounter{equation}{0} \setcounter{theorem}{0}

In \cite{Zla}, Zlato\v{s} proved, among several others, a result asserting
that, for a class of sparse random Sch\"{o}dinger operators $H_{\phi
}^{(\omega )}$on $\mathbb{Z}_{+}$, equivalent to the Jacobi matrices $%
J_{P,\phi }$ defined below, its spectral measure $\mu _{\phi }^{(\omega
)}(\lambda )$, restricted to the interval 
\begin{equation}
K=\left( -\sqrt{4-v^{2}/(\beta -1)},\sqrt{4-v^{2}/(\beta -1)}\right) ,\quad
v^{2}<4(\beta -1)  \label{K}
\end{equation}%
for a.e. disorder $\omega $ ($\omega =\{\omega _{n},~n\geq 1\}$ are
independent random variables uniformly distributed over $\left\{
-n,-n+1,\ldots ,n\right\} $) and a.e. phase boundary $\phi $, is purely
singular continuous (s.c.) with local Hausdorff dimension%
\begin{equation}
1-\frac{\ln \left( 1+\dfrac{v^{2}}{4-\lambda ^{2}}\right) }{\ln \beta }
\label{hd}
\end{equation}%
and dense pure point (p.p.) in the rest of the interval $\left[ -2,2\right] $%
. Although both the existence of p.p. spectrum and the sharpness of the
transition were mentioned without explicit proof (Theorem 6.3 of \cite{Zla}%
), such a proof exists \cite{R}. Theorem 6.3 is preceded by a two-page
argument establishing that the following statements hold for almost $\phi $
and for almost all points $\omega $ in the probability space: for almost all
energies $\lambda \in (-2,2)$ there is a subordinate generalised
eigenfunction with $O(n^{-\alpha _{\lambda }})$ decay and all other
generalised eigenfunctions have $n^{\alpha _{\lambda }}$ growth; $\alpha
_{\lambda }=\log \left( 1+v^{2}/\left( 4-\lambda ^{2}\right) \right) \left/
(2\log \beta )\right. $, and thus, for $\alpha _{\lambda }\leq 1/2$, i.e.,
the forthcoming interval $I$ of Theorem \ref{tmarsilwre}, assertion $(a)$
there holds, while in the complementary set\footnote{%
For $\alpha _{E}>1$ the subordinate generalized eigenfunction is an $l_{2}(%
\mathbb{Z}_{+})$ eingenfunction.} assertion $(b)$ there is valid (but
without the explicit characterization (\ref{sets})).

The excluded set $A$ of Lebesgue measure zero in $(-2,2)$, mentioned in the
above Theorem ($A$ has been written as the union of the two sets $A_{sc}$
and $A_{pp}$ situated in disjoint regions), is a common feature of both
methods of \cite{Zla} and ours, and is due to the peculiarities of the
definition of the essential support or minimal support of a measure (see
Definition 1 of \cite{GP}). Such set are not present in the cases of purely
p.p. or purely s.c. spectrum (see Theorems \ref{tsub1} and \ref{tult}) but
occurs here. The set $A_{sc}$ might be a discrete set or even dense p.p.,
the set $A_{pp}$ a countable union of a collection of (eventually different)
Cantor sets of zero Lebesgue measure, in which case both spectra would be
mixed. One distinctive feature of our method is that it allows us to make $%
A_{pp}$ explicit, which is seen also to be dense p.p., so that at least the
p.p. part of the spectrum is not mixed (see (\ref{sets})). We believe the
same happens for the s.c. part, but are not able to prove it.

Our objective is to prove both the existence of p.p. spectrum and the
sharpness of the transition with (s.c.) spectrum (Theorem \ref{tmarsilwre})
by a novel method, whose investigation started in \cite{MarWre} and
continued in \cite{SMW}. The existence proof involves the ergodicity of the
Pr\"{u}fer angles, corresponding to a solution $u$ of the eigenvalue
equation $J_{P,\phi}u=\lambda u$. As remarked by C. Remling \cite{Reml} in
his review of our Nonlinearity paper \cite{MarWre}, which introduced our
method, our new idea was to \textbf{fix} the energy and assume the Pr\"{u}%
fer angles $\left( \theta _{j}^{\omega }\right) _{j\geq 1}$ at $a_{j}$ are
uniformy distributed (u.d.) as a function of $j$ -- instead of the
traditional approach which exploits the u.d. of the Pr\"{u}fer angles in the
energy variable at fixed $a_{j}$. In \cite{MarWre} we were only able,
however, to fix the energy in the s.c. spectrum. In the present paper we are
able to bring our ideas to full fruition: we prove that the Pr\"{u}fer
angles are indeed u.d. (in the above sense) for any fixed energy, be it in
the s.c. or p.p. spectrum, by exploiting a (slight) modification of the
(optimal) metric extension of Weyl's criterion for uniform distribution by
Davenport, Erd\"{o}s and LeVeque \cite{DEL} (see also \cite{KN}, Theorem
4.2) -- this is our Theorem \ref{thethe}. Due to its simplicity and
naturalness, we expect that the method will find a wider range of
applications.

Our version of the sharpness of the spectral transition in Theorem 6.3 of 
\cite{Zla} proceeds along lines different from \cite{Zla}. Section \ref{ICPP}%
, devoted to the issue, replaces (\ref{EPP}) by a criterion for pure point
spectrum (Lemma \ref{l2}) that is suitable for sparse potentials. Lemma \ref%
{l2} is the analogue of Lemma 2.1 of \cite{Zla} in which the decay of a
subordinate solution is obtained.

It is worth noting that a sharp spectral transition occurs because both the
sparsity and the norm of transfer matrices grow exponentially with respect
to $j$. This conclusion is confirmed for models whose sparseness grow with
sub or super-exponential rates. For sub-exponential sparseness, the spectrum
is dense pure point (Theorem \ref{tsub1}); for super-exponential sparseness,
the spectrum is purely singular continuous (Theorem \ref{tult}). Varying the
intensity of perturbation, a model with super-exponential sparseness may
have purely singular continuous spectrum of Hausdorff dimension $0$ (Theorem %
\ref{tsub2}). It is, nevertheless, possible to choose $\delta =1$ in
relation (\ref{q}) in order to obtain a (possibly not sharp) transition from
singular continuous to pure point spectrum with the spectral measure being
of Hausdorff dimension $0$ (see Remark \ref{0st}).

We emphasize that no familiarity of the reader with \cite{MarWre} is
assumed. He (she) may look at Section $2$ for notation, Theorem \ref%
{tmarsilwre} for the statement of the main result and then concentrate on
Theorems \ref{teomet} and \ref{thethe}, Lemma \ref{l2} and the completion of
the proof of Theorem \ref{tmarsilwre} after Lemma \ref{l2}, which sketches
all basic steps of the proof. He (she) will turn eventually to pages $776$--$%
778$ of \cite{MarWre} to fill in some details of the derivation, which
follows the method of \cite{LS} and will not be repeated there.

The present paper is organized as follows. In Section \ref{SF} we present
the model, some definitions and results regarding its spectral transition.
In Section \ref{S2} we prove Theorem \ref{thethe} on the Pr\"{u}fer angles.
In Section \ref{ICPP} we improve the criterion for pure point spectrum.
Section \ref{SSES} contains Theorems \ref{tsub1}, \ref{tult} and \ref{tsub2}
which discuss subexponential and superexponential sparsity. In Section \ref%
{MO} we briefly relate the present result to a program started by Molchanov 
\cite{Mo} concerning an Anderson-like transition.

As a last remark, the Hausdorff dimension of the spectral measure is readily
obtained by applying the Gilbert--Pearson subordinacy criterion \cite{GP}
adapted for sparse potentials by Jitomirskaya and Last: see \cite{Zla} for
the diagonal case, \cite{SMW} for the nondiagonal case, by our method. We
refer to \cite{Reml} for additional references.

\section{Some Facts Regarding the Studied Class of Jacobi Matrices\label{SF}}

\setcounter{equation}{0} \setcounter{theorem}{0}

We begin with some definitions and basic results. The model analyzed here is
a small variation of the model studied in \cite{MarWre}, i.e., a class of
nondiagonal Jacobi matrices subject to exponentially sparse perturbations of
finite size. Although we consider this model for definiteness, because it is
part of a long-term project started in \cite{MarWre}, and further developped
in \cite{SMW}, we emphasize that the result apply with a few minor changes
to the diagonal model considered by Zlato\v{s} \cite{Zla} -- it is enough to
replace $(1-p^{2})/p$ in expression (\ref{r}) by the intensity parameter $v$
of equations (\ref{K}) and (\ref{hd}).

We consider off-diagonal Jacobi matrices 
\begin{equation}
J_{P}=\left( 
\begin{array}{ccccc}
0 & p_{0} & 0 & 0 & \cdots \\ 
p_{0} & 0 & p_{1} & 0 & \cdots \\ 
0 & p_{1} & 0 & p_{2} & \cdots \\ 
0 & 0 & p_{2} & 0 & \cdots \\ 
\vdots & \vdots & \vdots & \vdots & \ddots%
\end{array}%
\right) \;,  \label{MAJA}
\end{equation}%
\noindent for each sequence $P=(p_{n})_{n\geq 0}$ of the form%
\begin{equation}
p_{n}=\left\{ 
\begin{array}{lll}
p & \mathrm{if} & n=a_{j}^{\omega }\in \mathcal{A}\,, \\ 
1 & \mathrm{if} & \text{otherwise}\,,%
\end{array}%
\right.  \label{Pe}
\end{equation}%
\noindent for $p\in (0,1)$. The novelty with respect to \cite{MarWre}: $%
\mathcal{A}=\{a_{j}^{\omega }\}_{j\geq 1}$ is now a random set of natural
numbers $a_{j}^{\omega }=a_{j}+\omega _{j}$ with $a_{j}$ satisfying the
"sparseness" condition 
\begin{equation}
a_{j}-a_{j-1}=\beta ^{j}\;,\qquad \qquad j=2,3,\ldots  \label{spacon}
\end{equation}%
\noindent with $a_{1}+1=\beta \geq 2$ and $\omega _{j}$, $j\geq 1$,
independent random variables defined on a probability space $(\Omega ,%
\mathcal{B},\mu )$, uniformly distributed on the set $\Lambda
_{j}=\{-j,\ldots ,j\}$. These variables introduce uncertainty in the
position of the points where the $p_{n}$ differs from one. Such models are
nowadays called Poisson models, see \cite{SJ} and references given there.
Note that the support of $\omega _{j}$ increases only linearly with respect
to the index $j$.

The Jacobi matrices (\ref{MAJA}) can be written as 
\begin{equation*}
(J_{P}u)_{n}=p_{n}u_{n+1}+p_{n-1}u_{n-1}\;,
\end{equation*}%
\noindent with $u=(u_{n})_{n\geq 0}\in l_{2}(\mathbb{Z}_{+})$ and we denote
by $J_{P,\phi }$ the Jacobi matrix $J_{P}$ which satisfies a $\phi $%
-boundary condition at $-1$: 
\begin{equation}
u_{-1}\cos \phi -u_{0}\sin \phi =0\;.  \label{bouconphi}
\end{equation}

\begin{remark}
\label{r1} $J_{P,\phi }$ is a (noncompact) perturbation of the free Jacobi
matrix $J_{1,\phi }$, where $p_{n}=1$ for all $n\geq 0$: $J_{P,\phi
}=J_{1,\phi }+V_{P}$, the \textquotedblleft potential\textquotedblright\ $%
V_{P}$ composed by infinitely many random barriers whose distances grow
exponentially fast. A disordered potential of this sort was introduced by
Zlato\v{s} \cite{Zla}.
\end{remark}

\subsection{Essential spectrum}

Since, for some fixed realization of $\omega $, we are dealing with an
operator $J_{P,\phi }$ already studied on \cite{MarWre}, its essential
spectrum (for a definition, see section VII.3 of \cite{RESI}) remains the
same with probability $1$:

\begin{theorem}
\label{tee} The essential spectrum of the Jacobi matrix $J_{P,\phi }$,
defined by (\ref{MAJA}), (\ref{Pe}) and satisfying the $\phi $-boundary
condition (\ref{bouconphi}), for each realization of $\omega _{j}\in \Lambda
_{j}\equiv \{-j,\ldots ,j\}$, $j\geq 1$, is given by 
\begin{equation}
\sigma _{\text{ess}}=[-2,2]~.  \label{Ah}
\end{equation}
\end{theorem}

The proof of Theorem \ref{tee} is as the proof of Theorem 2.1 of \cite%
{MarWre}.

\subsection{Pr\"ufer variables}

To fix notation, we shall briefly review Sections 3 and 4 of \cite{MarWre}
incorporating eventual adjustments due to the randomness on the position of
the perturbations.

For $\lambda \in \mathbb{C}$, let 
\begin{equation}
T(n,n-1;\lambda )=\left( 
\begin{array}{cc}
\displaystyle\frac{\lambda }{p_{n}} & \displaystyle\frac{-p_{n-1}}{p_{n}}%
\vspace{2mm} \\ 
1 & 0%
\end{array}%
\right)  \label{mt}
\end{equation}%
\noindent be the $2\times 2$ transfer matrix associated to the $l_{2}(%
\mathbb{Z}_{+})$ solution of the eigenvalue equation $J_{P,\phi }u=\lambda u$%
. The equation 
\begin{equation}
\left( 
\begin{array}{c}
u_{n+1} \\ 
u_{n}%
\end{array}%
\right) =T(n,n-1;\lambda )\left( 
\begin{array}{c}
u_{n} \\ 
u_{n-1}%
\end{array}%
\right)  \label{uuTuu}
\end{equation}%
\noindent holds for $n\geq 0$, with ${\dbinom{u_{0}}{u_{-1}}}={\dbinom{\cos
\phi }{\sin \phi }}$ satisfying (\ref{bouconphi}) for some $\phi \in \lbrack
0,\pi ]$. The product of the $n+1$ first transfer matrices is denoted by 
\begin{equation}
T(n;\lambda )=T(n,n-1;\lambda )T(n-1,n-2;\lambda )\ldots T(0,-1;\lambda )\;.
\label{Tn}
\end{equation}

Given the definition (\ref{Pe}) of $p_{n}$ and the sparseness condition (\ref%
{spacon}), only three different $2\times 2$ matrices appear in the r.h.s. of
(\ref{Tn}):%
\begin{equation}
T_{-}(\lambda )=\left( 
\begin{array}{cc}
\lambda /p & -1/\vspace{2mm}p \\ 
1 & 0%
\end{array}%
\right) \ ,\qquad T_{+}(\lambda )=\left( 
\begin{array}{cc}
\lambda & -p\vspace{2mm} \\ 
1 & 0%
\end{array}%
\right) \ \mathrm{or}\qquad T_{0}(\lambda )=\left( 
\begin{array}{cc}
\lambda & -1\vspace{2mm} \\ 
1 & 0%
\end{array}%
\right)  \label{T-T+T0}
\end{equation}%
occurs depending on whether the left, the right or none of the two entries $%
n $ and $n-1$ in (\ref{mt}) are equal to $a_{j}^{\omega }\in \mathcal{A}$.
The free matrix $T_{0}(\lambda )$ is similar to a pure clockwise rotation by 
$\varphi \in (0,\pi )$: 
\begin{equation}
UT_{0}(\lambda )U^{-1}=\left( 
\begin{array}{cc}
\cos \varphi & \sin \varphi \vspace{2mm} \\ 
-\sin \varphi & \cos \varphi%
\end{array}%
\right) =R(\varphi )\;,  \label{R}
\end{equation}%
where $\lambda =2\cos \varphi $ and the similarity matrix is%
\begin{equation}
U\equiv \left( 
\begin{array}{cc}
0 & \sin \varphi \vspace{2mm} \\ 
1 & -\cos \varphi%
\end{array}%
\right) ~.\;  \label{U}
\end{equation}%
\noindent\ Note that $U$ is not uniquely defined since any other matrix $%
U^{\prime }=HU$, with $H$ commuting with $R$, satisfies (\ref{R}).

For every $n\in \mathbb{N}$ and $\omega _{j}\in \{-j,\ldots ,j\}$, $j\geq 1$%
, there is an integer $k$ such that $a_{k}^{\omega }\leq n<a_{k+1}^{\omega }$
and the conjugation of (\ref{Tn}) by $U^{-1}$, 
\begin{equation}
UT(n;\lambda )U^{-1}=R((n-a_{k}^{\omega })\varphi )P_{+-}(\lambda
)R((a_{k}^{\omega }-a_{k-1}^{\omega })\varphi )\cdots P_{+-}(\lambda
)R((a_{1}^{\omega })\varphi )\;,  \label{UTn}
\end{equation}%
\noindent intertwines $P_{+-}$ defined by 
\begin{equation*}
R(\varphi )P_{+-}(\lambda )R(\varphi )=UT_{+}(\lambda )T_{-}(\lambda
)U^{-1}\;
\end{equation*}%
with rotations.

Finally, (\ref{UTn}) are used together with (\ref{uuTuu}) and (\ref{Tn}) to
define the Pr\"{u}fer variables $(R_{k},\theta _{k}^{\omega })_{k\geq 0}$.
Observe that the Pr\"{u}fer angles are now random variables. The Pr\"{u}fer
radius, on the other hand, are random variables only as a function of the Pr%
\"{u}fer angles and their dependence on $\omega $ will be omitted. Like in
Section 3 of \cite{MarWre} (more specifically, equations (3.3)--(3.7)), if 
\begin{equation*}
v_{k}=(R_{k-1}\cos \theta _{k}^{\omega },R_{k-1}\sin \theta _{k}^{\omega })\;
\end{equation*}%
$k=1,2,\ldots $, are defined by 
\begin{equation*}
v_{k}=R((a_{k}^{\omega }-a_{k-1}^{\omega })\varphi )P_{+-}(\lambda )v_{k-1}
\end{equation*}%
\noindent for $k>1$ with $v_{1}=R(a_{1}^{\omega }\varphi )v_{0}$ where 
\begin{equation}
P_{+-}(\lambda )=\left( 
\begin{array}{cc}
p & 0\vspace{2mm} \\ 
\displaystyle\frac{1-p^{2}}{p}\cot \varphi & \displaystyle\frac{1}{p}%
\end{array}%
\right) ~,\;  \label{P+-}
\end{equation}%
then the Pr\"{u}fer variables satisfy a recursive relation (see equations (%
\ref{RN})-(\ref{PA}) below).

By equivalence of norms, the growth of $T(n;\lambda )$ may be controlled by
the euclidean norm (see Section 4 of \cite{MarWre} for detail): 
\begin{equation}
\Vert UT(n;\lambda )U^{-1}v_{0}\Vert ^{2}=\Vert UT(a_{N}^{\omega }+1;\lambda
)U^{-1}v_{0}\Vert ^{2}=R_{N}^{2}\;,  \label{R2N}
\end{equation}%
\noindent where the equality holds for any unit vector $v_{0}=(\cos \theta
_{0},\sin \theta _{0})$ and for each $n$ such that $a_{N}^{\omega }\leq
n<a_{N+1}^{\omega }$. Thus, from equations (\ref{R})-(\ref{P+-}), $R_{N}^{2}$
can be written as 
\begin{equation}
\left( R_{N}^{2}\right) ^{1/N}=\prod_{j=1}^{N}\left( \frac{R_{j}^{2}}{%
R_{j-1}^{2}}\right) ^{1/N}=\frac{1}{p^{2}}\exp \left\{ \frac{1}{N}%
\sum_{k=1}^{N}f(\theta _{k}^{\omega })\right\} \;,  \label{RN}
\end{equation}%
\noindent where 
\begin{equation}
f(\theta )=\ln \left( p^{4}\cos ^{2}\theta +(\sin \theta +(1-p^{2})\cot
\varphi \cos \theta )^{2}\right) \;  \label{ftheta}
\end{equation}%
\noindent and the Pr\"{u}fer angles $(\theta _{k}^{\omega })_{k\geq 1}$ are
defined recursively by 
\begin{equation}
\theta _{k}^{\omega }=\tan ^{-1}\left( \frac{1}{p^{2}}(\tan \theta
_{k-1}^{\omega }+\cot \varphi )-\cot \varphi \right) -(\beta ^{k}+\omega
_{k}-\omega _{k-1})\varphi  \label{PA}
\end{equation}%
\noindent for $k>1$ with $\theta _{1}^{\omega }$ given by 
\begin{equation*}
\theta _{1}^{\omega }=\theta _{0}-(a_{1}+\omega _{1})\varphi \;.
\end{equation*}

\subsection{Existence of singular continuous spectrum}

The nature of the spectrum of $J_{P,\phi }$ is intimately related to the
growth of the transfer matrix $T(n;\lambda )$. This connection was precisely
stated by Last and Simon \cite{LS}: the essential support $\Sigma _{\text{ac}%
}$ of the absolutely continuous part $\mu _{\text{ac}}$ of the spectral
measure $\mu =d\rho $ is given by 
\begin{equation}
\Sigma _{\text{ac}}=\left\{ \lambda :\liminf_{N\rightarrow \infty }\frac{1}{N%
}\sum_{n=0}^{N-1}\Vert T(n;\lambda )\Vert ^{2}<\infty \right\} \;.
\label{EAC}
\end{equation}%
\noindent By Theorem 1.6 of \cite{LS}, a sufficient condition for $\lambda
\in \sigma _{\text{sc}}$, the singular--continuous spectrum, is that 
\begin{equation}
\sum_{n=0}^{\infty }\Vert T(n;\lambda )\Vert ^{-2}=\infty \;,  \label{ESC}
\end{equation}%
\noindent since the eigenvalue equation $J_{P,\phi }u=\lambda u$ has no $%
l_{2}(\mathbb{Z}_{+})$--solutions. Theorem 1.7 of \cite{LS} asserts that the
eigenvalue equation has an $l_{2}(\mathbb{Z}_{+})$ solution (which is
precisely a sufficient condition for $\lambda \in \sigma _{\text{pp}}$, the
pure--point spectrum) if 
\begin{equation}
N_{\text{pp}}\equiv \sum_{n=0}^{\infty }\Vert T(n;\lambda )\Vert ^{2}\left(
\sum_{k=n}^{\infty }\Vert T(k;\lambda )\Vert ^{-2}\right) ^{2}<\infty \;.
\label{EPP}
\end{equation}

These relations were employed by Marchetti et. al. in \cite{MarWre} to
determine the spectral nature of the Jacobi matrices (with sparse
perturbations in deterministic positions). As we shall see, condition (\ref%
{EPP}) can be improved if Theorem 8.1 of \cite{LS}, used to establish (\ref%
{EPP}), is applied for the subsequence $T(a_{k}^{\omega }+1;\lambda )$, $%
k\geq 1$, instead of $T(n;\lambda )$ for any $n\in \mathbb{Z}_{+}$ (see
Section \ref{ICPP}). This improvement removes the gap in the spectrum
between the two intervals (\ref{I1}) and (\ref{I2}).

Provided the Pr\"{u}fer angles $(\theta _{k}^{\omega })_{k\geq 1}$ defined
by (\ref{PA}) are uniformly distributed modulo $\pi $, the Birkoff sum in
the r.h.s of (\ref{RN}) can be replaced, as $N$ tends to infinity, by an
integral over Lebesgue measure $d\theta $. Thus, the simple knowledge of the
integral of $f(\theta )$ gives us the asymptotic behavior of $R_{N}$.
Evidently, the assumption for this replacement has to be proven.

To formulate the questions addressed and enunciate our results, let us state
precisely what has been proved in \cite{MarWre}. We say that the Pr\"{u}fer
angles $(\theta _{k})_{k\geq 1}$ have uniform distribution modulo $\pi $
(u.d. mod $\pi $) if 
\begin{equation*}
\lim_{N\rightarrow \infty }\frac{\text{card}(\{k:\theta _{k}~\text{mod~}\pi
\in \lbrack 0,\theta ),1\leq k\leq N\})}{N}=\theta \;,
\end{equation*}%
\noindent holds for any $\theta \in \lbrack 0,\pi ]$, where $\text{card}(S)$
is the number of elements in $S$. By the ergodic theorem (see e.g. Theorem
1.1 of \cite{KN}), the sequence $(\theta _{k})_{k\geq 1}$ of real numbers is
u.d. mod $\pi $ if, and only if, 
\begin{equation}
\lim_{n\rightarrow \infty }\frac{1}{N}\sum_{k=1}^{N}f(\theta _{k})=\frac{1}{%
\pi }\int_{0}^{\pi }f(\theta )d\theta  \label{ergodic}
\end{equation}%
holds for every real-valued periodic Riemann integrable function $f$ of
period $\pi $.

It is easy to prove that $f$ defined by (\ref{ftheta}) satisfies all those
requirements and its integral is given by%
\begin{equation}
\frac{1}{\pi }\int_{0}^{\pi }f(\theta )d\theta =\ln \left( rp^{2}\right)
\label{lnrp2}
\end{equation}%
where%
\begin{equation}
r=1+\frac{(1-p^{2})^{2}}{4p^{2}}\text{cosec}^{2}\varphi ~.  \label{r}
\end{equation}%
(see e.g. equations (4.14)-(4.17) of \cite{MarWre}). Hence, the ergodic
theorem (\ref{ergodic}) combined with equation (\ref{RN}) and relations (\ref%
{EAC})-(\ref{EPP}), lead to Theorem 4.3 of \cite{MarWre}, which we state for
completeness:

\begin{theorem}
\label{tmarwre} Let $J_{P,\phi }$ be the Jacobi matrix given by (\ref{MAJA}%
), with the sequence $P=(p_{n})_{n\geq 0}$ defined by relations (\ref{Pe})
and (\ref{spacon}) and which satisfies the boundary condition (\ref%
{bouconphi}). Let 
\begin{equation}
I_{1}\equiv \left\{ \lambda \in \lbrack -2,2]:\frac{p^{2}}{\left(
1-p^{2}\right) ^{2}}(\beta -1)(4-\lambda ^{2})\geq 1\right\}  \label{I1}
\end{equation}%
\noindent 
\begin{equation}
I_{2}\equiv \left\{ \lambda \in \lbrack -2,2]:\frac{p^{2}}{\left(
1-p^{2}\right) ^{2}}(\beta ^{3}-1)(4-\lambda ^{2})<1\right\}  \label{I2}
\end{equation}%
\noindent with $p\in (0,1)$ and $\beta \in \mathbb{N}$, $\beta \geq 2$. If,
for each realization of $\omega $, the Pr\"{u}fer angles $(\theta
_{k}^{\omega })_{k\geq 1}$ are u.d. mod $\pi $ for $\varphi \in \lbrack
0,\pi ]\backslash A_{\theta _{0}}^{\omega }$, $A_{\theta _{0}}^{\omega }$ a
set of zero Lebesgue measure, possibly $\theta _{0},\omega $-dependent, then

(a) there exists a set $A_1$ of Lebesgue measure zero such that the spectrum
restricted to the set $I_1\backslash A_1$ is purely singular continuous,

(b) the spectrum of $J_{P,\phi }$ is pure point when restricted to $I_{2}$
for almost every $\phi \in \lbrack 0,\pi )$.
\end{theorem}

The proof of theorem \ref{tmarwre} is entirely contained in Section 4 of 
\cite{MarWre}.

For the deterministic model considered in \cite{MarWre} -- in particular,
for the model defined here for fixed $\omega $ -- only numerical evidence
for the uniform distribution of the Pr\"{u}fer angles $(\theta _{k}^{\omega
})_{k\geq 1}$ has been provided. The Pr\"{u}fer angles $(\theta _{k}^{\omega
})_{k\geq 1}$ may, however, be replaced by an u.d. mod $\pi $ sequence $%
(\zeta _{k})_{k\geq 1}$ where $\zeta _{k}(\varphi )$, as a function of $%
\varphi $, is a continuous piecewise linear interpolation of $\theta
_{k}(\varphi )$. It is shown in \cite{MarWre} that the error $E$ of
replacing their corresponding Birkoff sum can be made as small as one
wishes, provided the sparseness parameter $\beta >\beta _{0}(E,\lambda ,p)$
is large enough. The problem is that $(E,\lambda ,p)$ has to be fixed and $%
\beta >\beta _{0}$ makes the set $I_{2}$ given by (\ref{I2}) an empty set.
Theorem 5.8 of \cite{MarWre} states that only the singular continuous part
of Theorem \ref{tmarwre} is present for $\beta >\beta _{0}$ if the
hypothesis on the Pr\"{u}fer angles is removed.

\subsection{Transition from singular continuous to dense pure point spectrum}

In the present work we prove (Theorem \ref{thethe}) that the sequence of Pr%
\"{u}fer angles $(\theta _{k}^{\omega })_{k\geq 1}$ is u.d. mod $\pi $ for
all $\varphi \in \left[ 0,\pi \right] $ with exception of the set of
rational multiples of $\pi $ and for almost every $\omega $ with respect to
the product $\nu =\prod_{j\geq 1}\nu _{j}$ of uniform measures on $\Lambda
_{j}$. This, together with Theorem \ref{tmarwre}, establishes the existence
of a spectral transition from singular continuous to dense pure point
spectrum. We also prove in Section \ref{ICPP} a lemma under which Theorem %
\ref{tmarwre} holds with $\beta ^{3}$ in $I_{2}$ replaced by $\beta $. On
the set of parameters complementary to the set for which the singular
continuous spectrum occurs we have dense pure point spectrum. The spectral
transition is, consequently, sharp.

Our main result, the combination of Theorem \ref{thethe} and Lemma \ref{l2},
is summarized as follows:

\begin{theorem}
\label{tmarsilwre}Let $J_{P,\phi }$ be as in Theorem \ref{tmarwre}. Let $I$
be the interval (\ref{I1}) (for the diagonal case, replace $(1-p^{2})/p$ by $%
v$), $A=2\cos \pi \mathbb{Q}$ and set 
\begin{eqnarray}
A_{sc} &=&A\cap I  \notag \\
A_{pp} &=&A\cap I^{c}~  \label{sets}
\end{eqnarray}%
where $I^{c}$ means complementary of $I$ in $[-2,2]$. Then, for almost all $%
\omega $ with respect to uniform product measure on $\Lambda =\times
_{j=1}^{\infty }\left\{ -j,\ldots ,j\right\} $,

(a) the spectrum restricted to the set $I\backslash A_{sc}^{\prime }$, with $%
A_{sc}^{\prime }=A_{sc}\cup A^{\prime }$ and $A^{\prime }$ a set of Lebesgue
measure zero related with the definition of essential support of $\mu $, is
purely singular continuous;

(b) the spectrum of $J_{P,\phi }$ is dense pure point when restricted to $%
I^{c}\backslash A_{pp}$ for almost every $\phi \in \lbrack 0,\pi )$, where $%
\phi $ characterizes the boundary condition (\ref{bouconphi}). As we have
only excluded a countable set $A_{pp}$, the spectrum is purely p.p. in $%
[-2,2]\backslash I$.
\end{theorem}

\section{Uniform distribution of Pr\"{u}fer angles\label{S2}}

\setcounter{equation}{0} \setcounter{theorem}{0}

In this section we prove that the sequence of Pr\"{u}fer angles $(\theta
_{k}^{\omega })_{k\geq 1}$ is u.d. mod $\pi $ for almost every $\omega $ and
for all $\varphi \in \left[ 0,\pi \right] $ such that $\varphi /\pi $ is an
irrational number. For this purpose, we make use of the following slight
modification of the (optimal) metric extension of Weyl's criterion for
uniform distribution by Davenport, Erd\"{o}s and LeVeque \cite{DEL} (see
also Theorem 4.2 of \cite{KN}).

\begin{theorem}
\label{teomet}Let $(u_{n}(x))_{n\geq 1}$ be a sequence of real-valued random
variables defined in a probability space $\left( \Omega ,\mathcal{B},\mu
\right) $. For integers $h\neq 0$, $N\geq 1$ and $A\subset \Omega $, we set 
\begin{equation}
S_{h}(N,x)=\frac{1}{N}\sum_{n=1}^{N}e^{2ihu_{n}(x)}  \label{Sh}
\end{equation}%
\noindent and 
\begin{equation}
I_{h}(N,A)=\int_{A}|S_{h}(N,x)|^{2}d\mu (x)\;.  \label{IhN}
\end{equation}%
If the series $\sum_{N=1}^{\infty }I_{h}(N,A)/N$ converges for each $h\neq 0$%
, then the sequence $(u_{n}(x))$ is u.d. mod $\pi $ for almost all $x\in A$
with respect to $\mu $.
\end{theorem}

The random variables $\left\{ \omega _{j},\ j\geq 1\right\} $, defined on
probability space $\left( \Omega ,\mathcal{B},\mu \right) $, are
statistically independent and uniformly distributed in $\Lambda _{j}=\left\{
-j,\ldots ,j\right\} $. Let $\Lambda =\times _{j=1}^{\infty }\Lambda _{j}$
be the configuration space. The measure $\mu $ induces on the mesurable
space $(\Lambda ,\mathcal{F})$, where $\mathcal{F}$ is the natural product $%
\sigma $--algebra, a product measure 
\begin{equation*}
\nu (B)=\prod_{j=1}^{\infty }\nu _{j}(B_{j})=\mu \left( \omega
^{-1}(B)\right)
\end{equation*}%
where $\omega ^{-1}(B)=\displaystyle\bigcap_{j=1}^{\infty }\omega
_{j}^{-1}(B_{j})$, defined for all cylinder sets $B=\times _{j=1}^{\infty
}B_{j}\subset \Lambda $ with $\nu _{j}$ the uniform measure in $\Lambda _{j}$%
: $\nu _{j}(k)=1/(2j+1)$ for any $k\in \{-j,\ldots ,j\}$. We denote the
sequence of Pr\"{u}fer angles, defined by the recursive relations (\ref{PA}%
), either by $(\theta _{k}(x))_{k\geq 1}$ of by $(\theta _{k}^{\omega
})_{k\geq 1}$, depending on whether the probability space $\left( \Omega ,%
\mathcal{B},\mu \right) $ or $(\Lambda ,\mathcal{F},\nu )$ is referred.


The main result of this section is as follows.

\begin{theorem}
\label{thethe} The sequence of Pr\"{u}fer angles $(\theta _{k}(x))_{k\geq 1}$
($(\theta _{k}^{\omega })_{k\geq 1}$) is u.d. mod $\pi $ for all $\varphi
/\pi \in \lbrack 0,1]\backslash \mathbb{Q}$ and all $x\in \Omega $ ($\omega
\in \Lambda $) apart from a set with $\mu $ ($\nu $) measure $0$.
\end{theorem}

\noindent \textbf{Proof.} According to Theorem \ref{teomet}, we must show
that the series $\sum_{N=1}^{\infty }I_{h}(N,\Omega )/N$ converges, $%
I_{h}(N,\Omega )$ defined by (\ref{IhN}). It is, nevertheless, sufficient to
show that the series converges absolutely. Thus, by (\ref{PA}) and (\ref{Sh}%
), 
\begin{eqnarray}
I_{h}(N,\Omega ) &=&\int_{\Omega }|S_{h}(N,x)|^{2}d\mu (x)  \notag \\
&\leq &\frac{1}{N}+\frac{2}{N^{2}}\sum_{1\leq m<n\leq N}\left\vert
\int_{\Omega }e^{2ih(\theta _{m}(x)-\theta _{n}(x))}d\mu (x)\right\vert \;.
\label{esIh}
\end{eqnarray}%
\noindent Since $\theta _{m}^{\omega }$ with $m<n$ and $\widetilde{\theta }%
_{n}^{\omega }$, given by $\theta _{n}^{\omega }=g(\theta _{n-1}^{\omega
},\varphi )-(\beta ^{n}+\omega _{n}-\omega _{n-1})\varphi \equiv \tilde{%
\theta}_{n}^{\omega }-\omega _{n}\varphi $, are statistically independent of 
$\omega _{n}$, we have 
\begin{eqnarray}
\left\vert \int_{\Omega }e^{2ih(\theta _{m}(x)-\theta _{n}(x))}d\mu
(x)\right\vert &=&\left\vert \int_{\Lambda }e^{2ih(\theta _{m}^{\omega }-%
\tilde{\theta}_{n}^{\omega })}d\nu (\omega )\right\vert \left\vert
\int_{\Lambda _{n}}e^{2ih\omega _{n}\varphi }d\nu _{n}(\omega
_{n})\right\vert  \notag \\
&\leq &\left\vert \int_{\Lambda _{n}}e^{2ih\omega _{n}\varphi }d\nu
_{n}(\omega _{n})\right\vert .  \label{modint}
\end{eqnarray}%
\noindent The r.h.s. of (\ref{modint}) is the characteristic function of $%
\nu _{n}$ at $2h\varphi $ and 
\begin{eqnarray}
\left\vert \frac{1}{2n+1}\sum_{k=-n}^{n}e^{2ihk\varphi }\right\vert
&=&\left\vert \frac{1}{2n+1}\frac{\sin \left( 2n+1\right) h\varphi }{\sin
h\varphi }\right\vert  \notag \\
&\leq &\frac{1}{(2n+1)|\sin (h\varphi )|}<\infty \;,  \label{inwn}
\end{eqnarray}%
\noindent holds except if $\varphi $ is a rational multiple of $\pi $ in $%
\left[ 0,\pi \right] $, justifying the restriction stated in the theorem.

Together, equations (\ref{esIh}), (\ref{modint}) and (\ref{inwn}) yields 
\begin{equation*}
I_{h}(N,\Omega )\leq \frac{1}{N}+\frac{2}{|\sin h\varphi |N^{2}}\sum_{1\leq
m<n\leq N}\frac{1}{2n+1}<\frac{1}{N}\left( 1+\frac{1}{|\sin h\varphi |}%
\right) ,
\end{equation*}%
which implies that $\displaystyle\sum_{N=1}^{\infty }I_{h}(N,\Omega )/N$ is
finite for each $h\neq 0$, concluding the proof of Theorem \ref{thethe}.

\hfill $\Box $

\begin{remark}
If $\varphi /\pi $ is any rational number, there exists a $h\in \mathbb{Z}$, 
$h\neq 0$, such that $h\varphi /\pi =m\in \mathbb{Z}$ and we have 
\begin{equation*}
\lim_{\varphi ^{\prime }\rightarrow \varphi }\frac{\sin \left( 2n+1\right)
h\varphi ^{\prime }}{\sin h\varphi ^{\prime }}=2n+1\;,
\end{equation*}%
\noindent which implies that $I_{h}(N,\Omega )\leq \mathit{O}(1)$ and,
consequently, the estimate of $\sum_{N=1}^{\infty }I_{h}(N,\Omega )/N$
employed on the proof of Theorem \ref{thethe} diverges.
\end{remark}

\begin{remark}
\label{R3tt}The uniform assumption on $\nu _{j}$ is not necessary for our
analysis, but it is assumed for simplicity. To assure that the Pr\"{u}fer
angles $(\theta _{k}^{\omega })_{k\geq 1}$ are u.d. mod $\pi $ it is
sufficient to define the random variables $\omega _{j}$ supported in a
interval $\Lambda _{j}\equiv \lbrack -j^{\varepsilon },j^{\varepsilon }]\cap 
\mathbb{Z}$, $\varepsilon $ any positive real number. By the proof of
Theorem \ref{thethe} it is clear that in this case we would have 
\begin{equation*}
\sum_{N=1}^{\infty }\frac{I_{h}(N,\Omega )}{N}\leq \sum_{N=1}^{\infty }\frac{%
\mathit{O}(1)}{N^{1+\varepsilon }}<\infty \;,
\end{equation*}%
\noindent which, by Theorem \ref{teomet}, proves our assertion.
\end{remark}

\begin{remark}
\label{R4tt}The proof of Theorem \ref{thethe} independ on the sparseness
condition (\ref{spacon}). Although we are free to choose any function,
Theorems \ref{tmarwre} and \ref{tmarsilwre} require that 
\begin{equation}
a_{n}^{\omega }-a_{n-1}^{\omega }\geq 2  \label{2}
\end{equation}%
holds for all $n>1$ and $\omega \in \Lambda $.
\end{remark}

\section{Improved Criterion for Pure Point Spectrum\label{ICPP}}

\setcounter{equation}{0} \setcounter{theorem}{0}

The criterion (\ref{EPP}) for pure point spectrum is not optimal for sparse
Jacobi matrices. To obtain an interval stated in Theorem \ref{tmarsilwre} we
use instead the following

\begin{lemma}
\label{l2}Let $J_{P,\phi }$ be as in Theorem \ref{tmarwre} and let $%
t_{n}=\left\Vert T(a_{n}+1;\lambda )\right\Vert $ denote the spectral norm
of the associate transfer matrix $T(k;\lambda )$ at the point $k=a_{n}+1$
for some $\lambda \in \left[ -2,2\right] $. If%
\begin{equation}
\sum_{n=1}^{\infty }\beta ^{n}t_{n}^{-2}<\infty  \label{alpha}
\end{equation}%
and%
\begin{equation}
\sum_{n=1}^{\infty }\beta ^{n}t_{n}^{2}\left( \sum_{m=n}^{\infty
}t_{m}^{-2}\right) ^{2}<\infty  \label{beta}
\end{equation}%
are verified, then the eigenvalue equation $J_{P,\phi }u=\lambda u$ has an $%
l_{2}(\mathbb{Z}_{+})$ solution.
\end{lemma}

\noindent \textbf{Proof}\textit{.} We begin by adapting Theorem 8.1 of \cite%
{LS} for the model in consideration.

By (\ref{mt}) and (\ref{Pe}), 
\begin{equation*}
\left\Vert T(k,k-1;\lambda )\right\Vert ^{2}\leq \left\Vert T(k,k-1;\lambda
)\right\Vert _{E}^{2}\leq 1+\frac{1+\lambda ^{2}}{p^{2}}<\infty
\end{equation*}%
if $p\in (0,1)$, where $\left\Vert \cdot \right\Vert _{E}$ is the Euclidean
matrix norm, for $k,k-1\in \mathcal{A}$, otherwise $T_{j}(k,k-1;\lambda )$
is similar to a clockwise rotation $R(\varphi )$ by $\varphi =(1/2)\arccos
\lambda $ (see equation (\ref{R})). We write 
\begin{equation*}
T(a_{n}+1;\lambda )=A_{n}(\lambda )\cdots A_{1}(\lambda )
\end{equation*}%
where, for each $m\geq 2$ 
\begin{equation*}
A_{m}(\lambda )=T(a_{m}+1,a_{m};\lambda )\cdots
T(a_{m-1}+2,a_{m-1}+1;\lambda )=T_{-}T_{+}T_{0}^{\beta ^{m}-2}
\end{equation*}%
by (\ref{T-T+T0}). Denoting by $s_{n}=$ $\left\Vert A_{n}(\lambda
)\right\Vert $ the spectral norm of $A_{n}(\lambda )$, we have%
\begin{equation}
s_{n}\leq C\left( 1+\frac{1+\lambda ^{2}}{p^{2}}\right) \equiv B  \label{s}
\end{equation}%
$C=(1+\left\vert \cos \varphi _{j}\right\vert )/(1-\left\vert \cos \varphi
_{j}\right\vert )$, uniformly in $n$. As a consequence,%
\begin{equation}
\sum_{n=1}^{\infty }\frac{s_{n+1}^{2}}{t_{n}^{2}}<\infty  \label{sum}
\end{equation}%
verifies the assumption of Theorem 8.1 of \cite{LS} and provides the
existence of a subordinate solution $v$ associated to $\lambda $. Since the
sparseness parameter satisfies $\beta \geq 2$, (\ref{alpha}) is stronger
than (\ref{sum}) and that assumption is already verified under the
hypotheses of Lemma \ref{l2}.

The transfer matrices $T_{0}$, $T_{+-}:=T_{+}T_{-}$ given by (\ref{T-T+T0})
are, together with $T(a_{n}+1;\lambda )$ and $T^{\ast }(a_{n}+1;\lambda )$, $%
2\times 2$ unimodular real matrices. Hence the product $T^{\ast
}(a_{n}+1;\lambda )T(a_{n}+1;\lambda )$ is a $2\times 2$ unimodular
symmetric real matrix whose eigenvalues are $t_{n}^{2}$ and $t_{n}^{-2}$,
and corresponding (normalized) eigenvectors $\mathbf{v}_{n}^{+}$ and $%
\mathbf{v}_{n}^{-}$ are orthogonal: $\left( \mathbf{v}_{n}^{+},\mathbf{v}%
_{n}^{-}\right) =0$. We write $\mathbf{v}_{\phi }=\left( 
\begin{array}{c}
\cos \phi \\ 
\sin \phi%
\end{array}%
\right) $ and define $\phi _{n}$ by 
\begin{equation}
\mathbf{v}_{\phi _{n}}=\mathbf{v}_{n}^{-}~.  \label{vv}
\end{equation}%
Clearly, $\mathbf{v}_{n}^{+}=\mathbf{v}_{\phi _{n}+\pi /2}$ and by the
spectral theorem, we have%
\begin{eqnarray}
\left\Vert T(a_{n}+1;\lambda )\mathbf{v}_{\phi }\right\Vert ^{2} &=&\left( 
\mathbf{v}_{\phi },T^{\ast }(a_{n}+1;\lambda )T(a_{n}+1;\lambda )\mathbf{v}%
_{\phi }\right)  \notag \\
&=&t_{n}^{2}\left\vert (\mathbf{v}_{\phi },\mathbf{v}_{n}^{+})\right\vert
^{2}+t_{n}^{-2}\left\vert (\mathbf{v}_{\phi },\mathbf{v}_{n}^{-})\right\vert
^{2}  \notag \\
&=&t_{n}^{2}\sin ^{2}\left( \phi -\phi _{n}\right) +t_{n}^{-2}\cos
^{2}\left( \phi -\phi _{n}\right) ~.  \label{Ttt}
\end{eqnarray}%
$u=\left( u_{k}\right) _{k\geq 0}$, given by the second component $%
u_{k}=\left( T(k;\lambda )\mathbf{v}_{\phi }\right) _{2}$, solves the
eigenvalue equation $J_{P,\phi }u=\lambda u$ with boundary condition (\ref%
{bouconphi}).

Using properties of matrix norm together with (\ref{Ttt}) for $n+1$ and
definition (\ref{vv}), it can be shown (see proof of Theorem 8.1 of \cite{LS}%
)%
\begin{equation*}
\left\vert \phi _{n}-\phi _{n+1}\right\vert \leq \frac{\pi }{2}\frac{%
s_{n+1}^{2}}{t_{n}^{2}}~.
\end{equation*}

Condition (\ref{sum}) implies that the sequence $\left( \phi _{n}\right)
_{n\geq 1}$ has a limit $\phi ^{\ast }=\lim_{n\rightarrow \infty }\phi _{n}$%
. Hence, equation (\ref{Ttt}) and the telescope estimate%
\begin{equation*}
\left\vert \phi _{n}-\phi ^{\ast }\right\vert \leq \sum_{m=n}^{\infty
}\left\vert \phi _{m}-\phi _{m+1}\right\vert \leq \frac{\pi }{2}%
\sum_{m=n}^{\infty }\frac{s_{m+1}^{2}}{t_{m}^{2}}
\end{equation*}%
yields%
\begin{eqnarray}
\left\Vert T(a_{n}+1;\lambda )\mathbf{v}_{\phi ^{\ast }}\right\Vert ^{2}
&\leq &t_{n}^{2}\left( \phi ^{\ast }-\phi _{n}\right) ^{2}+t_{n}^{-2}  \notag
\\
&\leq &\frac{\pi }{2}B^{4}t_{n}^{2}\left( \sum_{m=n}^{\infty
}t_{m}^{-2}\right) ^{2}+t_{n}^{-2}  \label{tt}
\end{eqnarray}%
Note that, by definition (\ref{uuTuu}) of transfer matrix , $v_{k}=\left(
T(k;\lambda )\mathbf{v}_{\phi ^{\ast }}\right) _{2}$ is a strongly
subordinate solution in the sense that, for $u_{k}\equiv \left( T(k;\lambda )%
\mathbf{v}_{\phi ^{\ast }+\pi /2}\right) _{2}$, we have 
\begin{equation*}
\lim_{k\rightarrow \infty }\frac{\left\vert v_{k}^{2}+v_{k+1}^{2}\right\vert 
}{\left\vert u_{k}^{2}+u_{k+1}^{2}\right\vert }=0~,
\end{equation*}%
since $\left\Vert T(a_{n}+1;\lambda )\mathbf{v}_{\phi ^{\ast }+\pi
/2}\right\Vert ^{2}\geq t_{n}^{2}/2$ is satisfied for sufficiently large $n$
and, for some $n$ such that $a_{n}+1\leq k<a_{n+1}$, $\left\Vert T(k;\lambda
)\mathbf{v}_{\phi ^{\ast }}\right\Vert \leq B\left\Vert T(a_{n}+1;\lambda )%
\mathbf{v}_{\phi ^{\ast }}\right\Vert $ holds.

To conclude the proof, it remains to show that the subordinate solution is
also $l_{2}(\mathbb{Z}_{+})$ under the hypotheses (\ref{alpha}) and (\ref%
{beta}). By the equivalence of norms it is enough to show it for the $U$%
--norm $\left\Vert \cdot \right\Vert _{U}:=\left\Vert U\cdot \right\Vert $,
where $U$ is given by (\ref{U}) with $\lambda =2\cos \varphi $. For any $%
k\in \mathbb{N}$, we pick $n$ such that $a_{n}+1\leq k<a_{n+1}$ holds. We
have from (\ref{R2N}) 
\begin{equation*}
\left\Vert T(k;\lambda )\mathbf{v}_{\phi ^{\ast }}\right\Vert
_{U}^{2}=\left\Vert T(a_{n}+1;\lambda )\mathbf{v}_{\phi ^{\ast }}\right\Vert
_{U}^{2}\leq (1+\left\vert \cos \varphi \right\vert )\left\Vert
T(a_{n}+1;\lambda )\mathbf{v}_{\phi ^{\ast }}\right\Vert ^{2}
\end{equation*}%
which, by the sparseness condition (\ref{spacon}) together with (\ref{tt}),
yields 
\begin{equation*}
\sum_{k=1}^{\infty }\left\Vert T(k;\lambda )\mathbf{v}_{\phi ^{\ast
}}\right\Vert _{U}^{2}\leq B^{\prime }\sum_{n=1}^{\infty }\beta
^{n}t_{n}^{2}\left( \sum_{m=n}^{\infty }t_{m}^{-2}\right) ^{2}+B^{\prime
\prime }\sum_{n=1}^{\infty }\beta ^{n}t_{n}^{-2}
\end{equation*}%
for some constants $B^{\prime }$ and $B^{\prime \prime }$, concluding the
proof.

\hfill $\Box $

\noindent \textbf{Completion of the proof of Theorem \ref{tmarsilwre}.} By (%
\ref{RN})-(\ref{PA}), (\ref{ergodic}) and Theorem \ref{thethe}, we have for
all $\varphi /\pi \in \lbrack 0,1]\backslash \mathbb{Q}$ and a.e. $\omega $
(see the analogues of (\ref{ftheta}) and (\ref{PA}) for the diagonal case in 
\cite{Zla}), 
\begin{subequations}
\begin{equation}
\frac{1}{p^{2}}C_{N}^{-1/N}\exp \left( \int_{0}^{\pi }f(\theta )d\theta
\right) \leq \left( R_{N}^{2}\right) ^{1/N}\leq \frac{1}{p^{2}}%
C_{N}^{1/N}\exp \left( \int_{0}^{\pi }f(\theta )d\theta \right)  \label{eq.a}
\end{equation}%
where, by Koksma's inequalities \cite{K} (see also Theorem $5.1$ in Chap. $2$
of \cite{KN}),%
\begin{equation}
C_{N}\equiv \exp \left( -\frac{ND_{N}^{\ast }}{\pi }\int_{0}^{\pi
}\left\vert f^{\prime }(\theta )\right\vert d\theta \right) ~  \label{eq.b}
\end{equation}%
with $D_{N}^{\ast }$ the discrepancy of the sequence $\left( \theta
_{k}^{\omega }\right) _{k\geq 1}$ which tends to zero for u.d. sequences by
a theorem of Weyl \cite{W}) (see also Corolary $1.1$ in Chap. $2$ of \cite%
{KN} or \cite{DT}). This implies, by (\ref{eq.b}): 
\begin{equation}
\lim_{N\rightarrow \infty }C_{N}^{1/N}=1~.  \label{eq.c}
\end{equation}%
By (\ref{R2N}), (\ref{lnrp2}), (\ref{eq.a}) and the methods of \cite{LS}
(see \cite{MarWre} pp $776$--$778$ for a complete derivation) 
\begin{equation}
\tilde{C}_{n}^{-1}r^{n}\leq t_{n}^{2}\leq \tilde{C}_{n}r^{n}  \label{rtr}
\end{equation}%
holds for every $n$, with $r$ given by (\ref{r}) and $\tilde{C}_{n}$
satisfying (\ref{eq.c}). By (\ref{rtr}), (\ref{alpha}) and (\ref{beta}) are
simultaneously satisfied provided $\lambda \in \lbrack -2,2]\backslash 2\cos
\pi \mathbb{Q}$ and 
\end{subequations}
\begin{equation*}
\frac{\beta }{r}=\beta \left( 1+\frac{(1-p^{2})^{2}}{p^{2}(4-\lambda ^{2})}%
\right) ^{-1}<1~,
\end{equation*}%
which is the condition complementary to the one defining in (\ref{I1}).
Analogous results holds for the diagonal case, and this completes the proof
of Theorem \ref{tmarsilwre}.

\hfill $\Box $

\section{Sub and super--exponential sparseness\label{SSES}}

\setcounter{equation}{0} \setcounter{theorem}{0}

This section is devoted to the characterization of the spectral measure $\mu 
$ when the sparseness of the perturbations (\ref{Pe}) is sub or
super-exponential.

Our random set $\mathcal{A}=\{a_{n}^{\omega }\}_{n\geq 1}$ of natural
numbers $a_{n}^{\omega }=a_{n}+\omega _{n}$ is now admissible if 
\begin{equation}
a_{n}-a_{n-1}=\left[ e^{cn^{\gamma }}\right]  \label{gamma}
\end{equation}%
\noindent holds for every $n\geq 1$, with $c,\gamma >0$ and $\omega _{n}$, $%
n\geq 1$, independent random variables, uniform in $\Lambda _{n}$. Here $%
\left[ z\right] $ denotes the integer part of a real number $z$. The
sequence defined by (\ref{gamma}) increases sub or super--exponentially fast
depending on whether $\gamma <1$ or $\gamma >1$, and coincides with the one
previously defined if $\gamma =1$ and $c=\ln \beta $. Other sparseness
conditions with different growth rate may be dealt using the same methods
employed in this section but the family chosen is already wide enough for
our purposes.

Let the Jacobi matrix $J_{P,\phi }$ be defined by equations (\ref{MAJA}), (%
\ref{Pe}) and (\ref{bouconphi}) with the condition (\ref{spacon}) replaced
by (\ref{gamma}). We shall also allow $p$ in equation (\ref{Pe}) vary in
some cases. As the proof of Theorem \ref{tee} still holds in all these
situations (see proof of Theorem 2.1 of \cite{MarWre}), the essential
spectrum of $J_{P,\phi }$ remain the same as before: $\sigma _{\text{ess}%
}(J_{P,\phi })=\left[ -2,2\right] $.

\subsection{$\protect\gamma <1$ case}

We begin with the following

\begin{theorem}
\label{tsub1}Let $J_{P,\phi }$ be as in Theorem \ref{tmarwre} with the
sparseness condition (\ref{spacon}) replaced by (\ref{gamma}) with $0<\gamma
<1$ and $c>0$. Then 
\begin{equation}
\sigma _{\text{ess}}(J_{P,\phi })=\sigma _{\text{pp}}(J_{P,\phi })=\left[
-2,2\right] \;  \label{Ah1}
\end{equation}%
\noindent holds for almost every boundary condition $\phi \in \lbrack 0,\pi
) $ and almost all $\omega \in \Lambda = \times_{n=1}^{\infty }\{-j,\ldots
,j\} $ with respect to the product measure $\nu $.
\end{theorem}

\noindent \textbf{Proof.} This theorem is an extension of Theorem \ref%
{tmarwre} (b) for the case of sub-exponential sparseness. According to
Theorem 1.7 from \cite{LS}, we need to show that condition (\ref{EPP}) is
satisfied for every $\lambda \in \lbrack -2,2]$ and almost every boundary
condition $\phi \in \lbrack 0,\pi )$.\footnote{%
It is not necessary to apply the improved version Lemma \ref{l2}.} Theorem %
\ref{tee} then completes the proof.

If $n$ is such that $a_{N}^{\omega }\leq n<a_{N+1}^{\omega }$ is satisfied
for some $N\in \mathbb{N}$, by equations (4.19)-(4.23) of \cite{MarWre}, we
have that 
\begin{equation}
\sum_{k=n}^{\infty }\Vert T(k;\lambda )\Vert ^{-2}\leq
C^{2}\sum_{l=N}^{\infty }(e^{cl^{\gamma }}+2l)C_{l}r^{-l}\ \equiv S_{N}
\label{EPP1}
\end{equation}%
\noindent with $C=\sqrt{(1+|\cos \varphi |)/(1-|\cos \varphi |)}$, holds for 
$\lambda \in 2\cos ([0,\pi ]\backslash A^{\prime })$, $A^{\prime }$ a fixed
set of Lebesgue measure zero, and for all $\omega \in \Lambda $ apart a set
of $\nu $--measure zero. Note that Theorem \ref{thethe} assures that $%
C_{l}^{-1}r^{l}\leq R_{l}^{2}\leq C_{l}r^{l}$ holds with $r$ defined by (\ref%
{r}) and $C_{l}^{1/l}\searrow 1$ as $l\rightarrow \infty $.

Using (\ref{EPP}) together with (\ref{EPP1}) give 
\begin{equation*}
N_{\text{pp}}\leq C^{\prime }\sum_{N=1}^{\infty }(e^{cN^{\gamma
}}+2N)C_{N}^{+}r^{N}S_{N}^{2}\leq C^{\prime \prime }\sum_{N=1}^{\infty
}(e^{cN^{\gamma }}+2N)^{3}\left( C_{N}^{+}\right) ^{3}r^{-N}
\end{equation*}%
\noindent for some finite constants $C^{\prime }$ and $C^{\prime \prime }$.
Hence, $N_{\text{pp}}$ converges for all $\gamma $, $0<\gamma <1$, and for
every $\lambda $ in the set 
\begin{equation}
B=[-2,2]\backslash A_{1}\;,  \label{A1}
\end{equation}%
\noindent $A_{1}$ a set of Lebesgue measure zero, the eigenvalue equation $%
J_{P,\phi }u=\lambda u$ has an $l_{2}(\mathbb{Z}_{+})$ solution. Finally,
since (\ref{Ah}) remains true for every boundary condition $\phi \in \lbrack
0,\pi )$, Proposition 4.3 of \cite{MarWre} implies that $J_{P,\phi }$ has
only dense pure point spectrum in $[-2,2]$ for almost every $\phi \in
\lbrack 0,\pi )$. Thus, we have proved (\ref{Ah1}) and, consequently, the
theorem.

\hfill $\Box $

The pure point spectrum in Theorem \ref{tsub1} holds for any perturbation $%
0<p<1$ (see equation (\ref{Pe})). A less singular spectrum can be obtained
if we let $p_{a_{k}^{\omega }}$ tend to $1$ as $a_{k}^{\omega }$ tends to $%
\infty $.

\subsection{$\protect\gamma >1$ cases}

Let us consider now the Jacobi matrix $J_{P,\phi }$ subjected to the
super-exponential sparseness condition (\ref{gamma}) with $\gamma >1$. The
spectral measure of $J_{P,\phi }$ in this situation is given by Theorem
1.4(2) of \cite{Zla}.

\begin{theorem}
\label{tult}Let $J_{P,\phi }$ be defined as in Theorem \ref{tsub1} with $%
\gamma >1$ and $c>0$. Then the spectral measure $\rho $ of $J_{P,\phi }$ is
purely singular continuous and its Hausdorff dimension is 1 everywhere in $%
(-2,2)$ for almost all $\omega \in \Lambda = \times_{n=1}^{\infty
}\{-j,\ldots ,j\}$.
\end{theorem}

\noindent \textbf{Proof.} Together with the simple estimate%
\begin{equation}
\gamma (n+1)^{\gamma -1}\geq (n+1)^{\gamma }-n^{\gamma
}=\int_{n}^{n+1}\gamma x^{\gamma -1}dx\geq \gamma n^{\gamma -1}~,
\label{ineq}
\end{equation}%
equations (\ref{Pe}) and (\ref{gamma}), yields 
\begin{equation*}
\frac{a_{n}^{\omega }}{a_{n+1}^{\omega }}\leq \frac{e^{cn^{\gamma }}+2n}{%
e^{c(n+1)^{\gamma }}-2n}=\frac{1+2ne^{-cn^{\gamma }}}{e^{c\gamma n^{\gamma
-1}}-2ne^{-cn^{\gamma }}}\longrightarrow 0
\end{equation*}%
\noindent as $n$ tends to infinity, for all $\omega \in \Lambda $, $\gamma
>1 $ and $c>0$. Thus, by Theorem 1.4(2) of \cite{Zla}, the Hausdorff
dimension of $\rho $ is $1$.

\hfill $\Box $

\begin{remark}
The techniques developed in \cite{MarWre} can also be used to prove that $%
\rho $ has Hausdorff dimension 1 (see \cite{SMW}). \noindent
\end{remark}

The sparse perturbation may be chosen so that the spectral measure, despite
of being singular continuous, has Hausdorff dimension 0. For this, let the
sequence $P=(p_{n})_{n\geq 0}$ vary depending on the $a_{k}^{\omega }\in 
\mathcal{A}$: 
\begin{equation}
p_{n}=\left\{ 
\begin{array}{lll}
q_{k} & \mathrm{if} & n=a_{k}^{\omega }\in \mathcal{A}\,, \\ 
1 & \mathrm{if} & \text{otherwise}\,%
\end{array}%
\right. ~  \label{Pe1}
\end{equation}%
where $q_{k}$ goes to zero, as $k\rightarrow \infty $, super--exponentially
fast: 
\begin{equation}
e^{cn^{\gamma }-c_{1}n^{\delta }}\leq \prod_{k=0}^{n}q_{k}^{-2}\leq
e^{cn^{\gamma }-c_{2}n^{\delta }}~  \label{q}
\end{equation}%
holds for some $c_{1}\geq c_{2}>0$ and $\gamma >\delta >1$ such that $\delta
>\gamma -1$ holds. In this case, $f$ defined by (\ref{ftheta}) is no longer
the same function in (\ref{RN}) and, for each $k=1,2,\ldots $, we have%
\begin{eqnarray*}
f_{k}(\theta ) &=&\ln \left( q_{k}^{4}\cos ^{2}\theta +(\sin \theta
+(1-q_{k}^{2})\cot \varphi \cos \theta )^{2}\right) \\
&\equiv &\ln \left( a_{k}+b_{k}\cos 2\theta +c_{k}\sin 2\theta \right)
\end{eqnarray*}%
where%
\begin{eqnarray*}
2a_{k} &=&\left( 1-q_{k}^{2}\right) ^{2}\cot ^{2}\varphi +1+q_{k}^{4} \\
2b_{k} &=&\left( 1-q_{k}^{2}\right) ^{2}\cot ^{2}\varphi -1+q_{k}^{4} \\
c_{k} &=&\left( 1-q_{k}^{2}\right) \cot \varphi
\end{eqnarray*}%
satisfies%
\begin{equation*}
a_{k}^{2}=b_{k}^{2}+c_{k}^{2}+q_{k}^{4}~.
\end{equation*}%
An explicit calculation, yields%
\begin{equation}
a_{k}-\sqrt{a_{k}^{2}-q_{k}^{4}}\leq a_{k}+b_{k}\cos 2\theta +c_{k}\sin
2\theta \leq a_{k}+\sqrt{a_{k}^{2}-q_{k}^{4}}=\sin ^{-2}\varphi
(1+O(q_{k}^{2}))  \label{qq}
\end{equation}%
which can be used to draw the following conclusion

\begin{theorem}
\label{tsub2}Let $J_{P,\phi }$ be defined by (\ref{MAJA}), (\ref{Pe1}), (\ref%
{gamma}) and (\ref{bouconphi}) with $\gamma >1$ and $c>0$. Then, for almost
every boundary condition $\phi \in \lbrack 0,\pi )$ and all $\omega \in
\Lambda = \times_{n=1}^{\infty }\{-j,\ldots ,j\}$, the spectral measure $%
\rho $ of $J_{P,\phi }$ restricted to $(-2,2)$ is purely singular continuous
and has $0$ Hausdorff dimension.
\end{theorem}

\noindent \textbf{Proof.} Let us begin proving that the essential spectrum
of $J_{P,\phi }$ is entirely singular continuous. According to Theorem 1.6
from \cite{LS}, it is sufficient to show that (\ref{ESC}) holds for every $%
\lambda \in \lbrack -2,2]$ apart a set of zero Lebesgue measure. From (\ref%
{R2N}), (\ref{RN}), (\ref{q}) and (\ref{qq}),

\begin{equation}
\sum_{n=0}^{\infty }\Vert T(n;\lambda )\Vert ^{-2}\geq \tilde{C}%
^{-2}\sum_{n=0}^{\infty }e^{cn^{\gamma }}\prod_{k=0}^{n}q_{k}^{2}\sin
^{2n}\varphi \geq \tilde{C}^{-2}\sum_{n=0}^{\infty }e^{c_{2}n^{\delta }}\sin
^{2n}\varphi  \label{Tsum}
\end{equation}%
diverges for $\varphi \neq 0,\pi $, which establishes the claim.

Now we apply Corollaries 4.2 and 4.5 of \cite{JitLast}. Suppose that for
some $\alpha \in \lbrack 0,1)$ and every $\lambda $ in some Borel set $A$, 
\begin{equation}
\limsup_{l\rightarrow \infty }\frac{1}{l^{2-\alpha }}\sum_{n=0}^{l}\left%
\Vert T(n;\lambda )\right\Vert ^{2}<\infty \;.  \label{dnor}
\end{equation}%
\noindent Then the restriction $\rho (A\cap \cdot )$ of the spectral measure 
$\mu =d\rho $ is $\alpha $--continuous. Suppose, on the other hand, that 
\begin{equation}
\liminf_{l\rightarrow \infty }\frac{\left\Vert u_{\mathrm{sub}}\right\Vert
_{l}^{2}}{l^{\alpha }}=0  \label{c4.5a}
\end{equation}%
\noindent for every $\lambda $ in some Borel set $A$, where $u_{\mathrm{sub}%
} $ is a subordinate solution of $J_{P,\phi }u=\lambda u$. Then the
restriction $\rho (A\cap \cdot )$ is $\alpha $--singular.

By a calculation analogous to (\ref{Tsum}), for $a_{m}<l\leq a_{m+1}$ 
\begin{eqnarray*}
\sum_{n=0}^{l}\Vert T(n;\lambda )\Vert ^{2} &\leq
&C^{2}\sum_{k=1}^{m}e^{ck^{\gamma }}\prod_{j=0}^{k}q_{j}^{-2}\sin
^{-2k}\varphi +C^{2}(l-a_{m})\prod_{j=0}^{m}q_{j}^{-2}\sin ^{-2m}\varphi \\
&\leq &C^{2}\sum_{k=1}^{m}e^{2ck^{\gamma }}e^{-c_{2}k^{\delta }}\sin
^{-2k}\varphi +C^{2}(l-a_{m})e^{cm^{\gamma }}e^{-c_{2}m^{\delta }}\sin
^{-2m}\varphi \\
&\leq &C^{\prime }\left( \sum_{k=0}^{m}e^{2ck^{\gamma }}+\left(
l-a_{m}\right) e^{cm^{\gamma }}\right) \leq C^{\prime \prime }l^{2}
\end{eqnarray*}%
implies that $\rho $ is $\alpha $--continuous only for $\alpha =0$. Since,
by Lemma 2.1 of \cite{Zla}, a subordinate solution of $J_{P,\phi }u=\lambda
u $, $u_{\mathrm{sub}}(n)$, decays as fast as the transfer matrix $\Vert
T(n;\lambda )\Vert $ grows, we have%
\begin{equation}
\left\Vert u_{\mathrm{sub}}\right\Vert _{l}^{2}\leq
C^{-2}\sum_{n=0}^{l}\Vert T(n;\lambda )\Vert ^{-2}\leq
C^{-2}\sum_{k=0}^{m+1}e^{c_{1}k^{\delta }}\sin ^{2k}\varphi \leq C^{\prime
}e^{c_{1}m^{\delta }}~  \label{usub}
\end{equation}%
for some finite constant $C^{\prime }$. Since $l^{\alpha }\geq e^{c\alpha
m^{\gamma }}$ and $\gamma >\delta $, the inferior limit (\ref{c4.5a}) goes
to $0$ for any $\alpha >0$ and the spectral measure $\rho $ is $\alpha $%
-singular in the whole interval $\left( -2,2\right) $ for any $\alpha $
concluding the proof of the theorem.

\hfill $\Box $

\begin{remark}
Note that $q_{k}^{-2}=e^{c\gamma k^{\gamma -1}-c_{1}\delta k^{\delta -1}}$,
with $\gamma >\delta >1$, $\delta >\gamma -1$ and $c_{1}$ sufficiently
large, satisfies (\ref{q}), in view of%
\begin{equation*}
n^{\kappa }=\kappa \int_{0}^{n}x^{\kappa -1}dx<\sum_{k=0}^{n}\kappa
k^{\kappa -1}<\kappa \int_{0}^{n+1}x^{\kappa -1}dx=(n+1)^{\kappa }
\end{equation*}%
for $\kappa =\gamma $ and $\delta $ together with inequality (\ref{ineq}).
\end{remark}

\begin{remark}
\label{0st}It follows from equations (\ref{Tsum}) and (\ref{usub}) with $%
\delta =1$ that the spectral measure restricted to $\left\{ \lambda \in %
\left[ -2,2\right] :\left( 1-\lambda ^{2}/4\right) e^{c_{2}}>1\right\} $ is
singular continuous with $0$ Hausdorff dimension and pure point when
restricted to $\left\{ \lambda \in \left[ -2,2\right] :\left( 1-\lambda
^{2}/4\right) e^{c_{1}}<1\right\} $.
\end{remark}

\section{Conclusions and outlook\label{MO}}

The transition depicted in Theorem \ref{tmarsilwre} may be considered as an
Anderson-like transition in one dimension, along the lines of the
(multidimensional) program laid out by Molchanov \cite{Mo}. Of particular
relevance in this context is the important fact, shown by \cite{DJLS} that
Anderson localization is unstable under rank one perturbations, i.e., a rank
one perturbation may change the p.p. spectrum into s.c..The same authors
point out, however, that this s.c. spectrum must have zero Hausdorff
dimension. Since the (local) Hausdorff dimension of the s.c. spectrum in
Theorem \ref{tmarsilwre} is (a.e.) nonzero, we conclude that the transition
is robust. Also in this connection, since we are close to the difference
Laplacian (in this regime of strong sparsity), what is really surprising in
Theorem \ref{tmarsilwre} is the existence of the pure point spectrum. One
interesting possibility of pursuing this program is to analyse the Kronecker
sum of two or more copies of $J_{P,\phi }$ \cite{SMW1}. Following an idea of
Malozemov-Molchanov \cite{MaMo}, the continuous part of the spectral measure
may turn out to be absolutely continuous: this has been established,
however, only for some classes of potentials with superexponential sparsity (%
\cite{S},\cite{KR}, \cite{SMW1}), for which, however, the spectrum is purely
s.c. (see Theorem \ref{tult}), and proving that the continuous part of
models with mixed spectrum in two or more dimensions, such as the Kronecker
sum of two or more copies of the present or similar models, becomes
absolutely continuous, presents a challenging open problem.

\end{document}